\documentclass[12pt]{article}

\usepackage{amsmath}
\usepackage{amssymb}
\usepackage{amscd}
\usepackage{latexsym}
\usepackage{theorem}
\usepackage{setspace}
\usepackage{epsfig}
\usepackage{psfrag}
\usepackage{amsfonts}
\usepackage{enumerate}
\usepackage{euscript}

\newtheorem{theorem}{Theorem}[section]
\newtheorem{lemma}[theorem]{Lemma}
\newtheorem{proposition}[theorem]{Proposition}
\newtheorem{corollary}[theorem]{Corollary}

{\theorembodyfont{\rmfamily} \theoremstyle{plain}
\newtheorem{defn}[theorem]{Definition}

\newtheorem{remark}[theorem]{Remark}

}

\numberwithin{equation}{section} \numberwithin{figure}{section}

\newcommand{\Z}{\mathbb{Z}}
\newcommand{\CP}{\mathbb{CP}}
\newcommand{\C}{\mathbb{C}}
\newcommand{\R}{\mathbb{R}}

\newcommand{\spa}[1]{\ensuremath{ P^{#1}( X, \beta,J)^{*}}}
\newcommand{\tpa}[1]{\ensuremath{ W^{#1}( f^{*}TX)}}

\newcommand{\am}[2]{\ensuremath{ M^{#1}_{#2}(X,\beta,J)}}
\newcommand{\sam}[2]{\ensuremath{ M^{#1}_{#2}(X,\beta,J)^{*}}}
\newcommand{\cam}[2]{\ensuremath{ \overline{M}^{#1}_{#2}(X,\beta,J)}}
\newcommand{\an}{\ensuremath{\Lambda^{0,1} }}
\newcommand{\zdc}[2]{\ensuremath{H_{D}^{0}(#1, #2)}}
\newcommand{\fdc}[2]{\ensuremath{H_{D}^{1}(#1, #2)}}
\newcommand{\red}{\ensuremath{\beta_{1} \overset{r,s}{\bigvee}
\beta_{2}}}

\begin{document}
\noindent {\LARGE \bf Intersection theoretic properties on the\\
\medskip moduli space of genus 0 stable maps to
\\ \medskip a semipositive symplectic 4-manifold }

\vspace{1cm}

\begin{center}
 {\bf Seongchun Kwon }\footnote{\mbox{~~~~}{\it MSC 2000
Subject Classification}: Primary: 53D30; 53D45,  \hspace{0.1in}
Secondary: 14C17 \newline \mbox{~~~~}{\it Keywords}: transversality, Gromov-Witten invariant, Moduli space of stable maps}\\
\end{center}

\begin{abstract}
We characterize transversality, non-transversality properties on the
moduli space of genus 0 stable maps to a semipositive symplectic
manifold of dimension 4, when GW([point], \ldots,[point]) is
enumerative. In particular, we show that the intersection theoretic
property depends on the existence of a critical point on a stable
map.
\end{abstract}

\section{Introduction}\label{s;intro}

The Gromov-Witten invariant is defined as an integration over a
moduli space. Let $[a_{i}]$ be a cohomology class which is
Poincar\`{e} dual to a point $a_{i}$, $i=1, \ldots,
c_{1}(f^{*}TX)-1$, in a compact semipositive symplectic manifold $X$
of dimension 4. We call the Gromov-Witten invariant GW($[a_{1}],
\ldots, [a_{n}]$) enumerative if it is a positive integer and counts
the number of stable maps passing through $c_{1}(f^{*}TX)-1$ points
in general position. It implies the following:
\begin{itemize}
\item $ev_{i}^{-1}(a_{i})$, $i=1, \ldots, n:=c_{1}(f^{*}TX)-1$, meets transversally.
\item The
number of points in $\bigcap_{i=1}^{n} ev_{i}^{-1}(a_{i})$ doesn't
vary depending on the general choice of configuration points $a_{1},
\ldots, a_{n}$.
\end{itemize}

\par In 2003, during his invitation, Gang Tian predicted the intersection theoretic properties on the
moduli space of stable maps when the target space is $\CP^{2}$. His
conjecture relates the non-transversality properties of the cycles
$ev_{i}^{-1}(a_{i})$, $i=1, \ldots, n$, with the properties of the
stable maps which represent the points in $\bigcap_{i=1}^{n}
ev_{i}^{-1}(a_{i})$. The main results in this paper prove his
conjecture.\\
\indent Usually, it is very hard to calculate the intersection
multiplicity at an intersection point in $\bigcap_{i=1}^{n}
ev_{i}^{-1}(s_{i})$ straightforwardly. In this paper, we calculate
the intersection multiplicities in a tricky way. Some technical
motivations came from Tian's suggestion. He related calculations to
the deformation properties of stable maps. That is, an appearance of
a skyscraper sheaf at a critical point can be used to show the
intersection cycles' non-transversality properties. Practically,
studies on the local structure of a moduli space with a fixed
generic almost complex structure $J$ and the singularity analysis of
the product of the $i$-th evaluation maps enabled the author to
prove his
conjecture. \\

\vspace{2mm}

 The main results of this paper are:

\begin{theorem}\label{t;maintransversality}
Let $n := c_{1}(f^{*}TX) -1$. Let $\mathbf{f}$ be in
$\bigcap_{i=1}^{n} ev_{i}^{-1}(q_{i})$, where $q_{i}$,
$i=1, \ldots, n$, are points in general position in the compact
semipositive symplectic 4-manifold $X$. Then, the following holds.\\
(i) If $\mathbf{f}$ is represented by a stable map which is an
immersion and has an irreducible domain curve,
 then the intersection multiplicity at $\mathbf{f}$ is one.\\
(ii)If $\mathbf{f}$ is represented by a cuspidal stable
map(Definition in \ref{d;maps}) whose marked points are not critical
points, then the intersection multiplicity at $\mathbf{f}$ is two.
\end{theorem}

\begin{theorem}\label{t;cusmul} Let $n := c_{1}(f^{*}TX) -1$.
Let $X$ be a compact semipositive symplectic 4-manifold. The
cuspidal stable maps locus is the unique equi-singular locus in
$\am{}{n}$ of real codimension $\leq 2$ on which transversality
uniformly fails.
\end{theorem}

\par This paper aims to exhibit the symplectic counterpart of the paper
\cite{kwona} in algebraic geometry category. For that purpose, we
express the local structure of the moduli space algebraically.
Preparations are done in sec.\ref{s;d-cohomology} and
sec.\ref{ss;normalchom}. The author doesn't claim any new results in
sec.\ref{s;d-cohomology} and sec.\ref{ss;normalchom}. This paper's
main part starts from sec.\ref{ss;tansplit}. In
sec.\ref{ss;tansplit}, we calculate the structure of the tangent
space of $\cam{}{n}$. In sec.\ref{s;index}, we calculate the
singularities of the $ev$ map and the index at the singularities.
Some results in this section look similar to some results in
\cite{she}. However, the moduli space we consider in this paper is
different. Shevchishin considered the total moduli space and the
singularity analysis were done on the total moduli space. If we fix
the almost complex structure $J$ in the total moduli space in
\cite{she} as in this paper, then the moduli space Shevchishin
worked become a set of discrete points where the singularity
analysis is not possible. In Sec.\ref{s;transversalityofpoint}, we
prove Tian's conjecture on the transversality properties of the
cycles $ev_{i}^{-1}(q_{i})$, $i = 1, \ldots, c_{1}(f^{*}TX) -1$.

\section{Shevchishin's D-cohomology group}\label{s;d-cohomology}

Let $(X, \omega, J)$ be a compact semipositive symplectic
2$n$-dimensional manifold with a $\omega$-compatible almost complex
structure $J$ with $C^{l}$-smooth,$l \geq 2$. Assume $kp >2$ and $k
\leq l$. \emph{A parameter space of simple $J$-holomorphic maps},
denoted by $\spa{k,p}$, is a space of continuous maps $f: \CP^{1}
\rightarrow X$ such that the $k$-th derivative of $f$ is of class
$L^{p}$ and
$f_{*}([\CP^1])= \beta \in H_{2}(X;\Z)$. \\

Let $L^{p}(\CP^1, \an f^{*}TX)$ be a Banach space of
$L^{p}$-integrable $f^{*}TX$-valued $(0,1)$-forms on $\CP^1$. Then,
$\bigsqcup_{f} L^{p}(\CP^1, \an f^{*}TX)$ is an infinite dimensional
vector bundle over the space $P^{1,p}(X, \beta)$ of continuous maps
$f: \CP^{1} \rightarrow X$ of a class $L^{1,p}$ representing the
homology class $\beta$. Let $\overline{\partial}_{J}$ be a complex
anti-linear section on $P^{1,p}(X, \beta)$:
\[ \overline{\partial}_{J}(f) := \frac{1}{2}(df + J \circ df \circ
j) \in L^{p}(\CP^1, \an  f^{*}TX)  \]

Then, the subspace $\spa{1,p}$ is an open subset of the
intersection with the zero section $\bar{\partial}_{J}^{-1}(0)$.\\

\vspace{2mm}

The linearized operator $D_{f,J}$ of $\overline{\partial}_{J}$ is:

\begin{gather}
D_{f,J}: \tpa{1,p} \rightarrow  L^{p}(\CP^1, \an
f^{*}TX)   \notag \\
D_{f,J}(v) = \frac{1}{2}(\nabla v + J \circ \nabla v \circ j +
\nabla_{v} J \circ df \circ j ). \label{eq;D}
\end{gather}

\vspace{2mm}

$D_{f,J}$ is an elliptic first order partial differential operator.
The elliptic operator $D_{f,J}$ defines a two-step elliptic complex:

\begin{equation}
 0 \rightarrow \tpa{1,p} \stackrel{D_{f,J}}{\rightarrow}
L^{p}( \an f^{*}TX)   \rightarrow 0. \label{eq;ell}
\end{equation}

The above two-step elliptic complex defines cohomology groups, which
were named as \emph{$D$-cohomology groups} by Shevchishin.

\[ H^{0}_{D}(S, f^{*}TX) := \mbox{Ker}D_{f,J}, \hspace{1cm} H^{1}_{D}(S, f^{*}TX) := \mbox{Coker}D_{f,J},      \]

where $S$ is a Riemann surface. $H^{0}_{D}(S, f^{*}TX)$ and $
H^{1}_{D}(S, f^{*}TX)$ are finite dimensional vector spaces over
$\R$ because $D_{f,J}$ is a Fredholm operator. We will call an
element $v$ in
$H^{0}_{D}(S, f^{*}TX)$ a \emph{pseudo-holomorphic section}.\\

\begin{remark}\label{r;indep} $D$-cohomology groups are defined with any
$k,p$ if $kp >2 $ and $k \leq l$, where $J$ is $C^l$-smooth. The
elliptic regularity implies that \zdc{\CP^1}{f^{*}TX},
\fdc{\CP^{1}}{f^{*}TX} don't depend on the choice of the functional
space. So, the resulting $D$-cohomology groups are independent of
$k,p$.
\end{remark}

For generic $J$, the linearized operator $D_{f,J}$ at $f$ is
surjective if $f$ is a simple map. It implies that the 1st
D-cohomology group vanishes. The following Proposition is from the
Implicit function theorem. See \cite[Theorem A.3.3]{mas} or
\cite[Lemma 2.2.4]{she}.

\begin{proposition}\label{c;tanofpa}
$\spa{1,p}$ is a smooth separable Banach manifold whose tangent
space at $f$ is:
\[T_{f} \spa{1,p} = H^{0}_{D}(\CP^1, f^{*}TX ) \]

\end{proposition}

\section[Local Structure]{Local structures of the moduli space of
genus 0 stable maps}\label{s;local}

\subsection{D-cohomology Groups for Normal
Sheaf}\label{ss;normalchom}

\begin{lemma}
Let $X$ be a compact semipositive symplectic 4-manifold. Let $f$ be
a $J$-holomorphic map such that $f_{*}([\CP^{1}])$ is non-trivial.
Let $\mathcal{O}(V)$ be a sheaf of sections of a vector bundle $V$.
Then, the sequence of coherent sheaves

\begin{equation}\label{e;short}
0 \rightarrow \mathcal{O}(T \CP^1) \stackrel{df}{\rightarrow}
\mathcal{O}(f^{*}TX) \rightarrow
\mathcal{O}(f^{*}TX)/df(\mathcal{O}(T \CP^1)) \rightarrow 0
\end{equation}

is exact.
\end{lemma}

Proof. It is obvious that the sequence is surjective to
$\mathcal{O}(f^{*}TX)/df(\mathcal{O}(T \CP^1))$ and exact at
$\mathcal{O}(f^{*}TX)$. Thus, it is enough to show that $df$ is
injective. Since $f$ is a $J$-holomorphic map and $J$ is in the
class $C^{l}$, $l \geq 2$, the number of critical points of $f$ is
finite. See \cite{maw}. Thus, $f$ is not locally constant by the
unique continuation theorem. It implies that $df(v) = 0$ only if $v$
is a trivial section.   \hfill
$\Box$\\

\vspace{5mm}

Let $c \in \CP^1$ be a critical point of the $J$-holomorphic map
$f$. Let $O(c)$ be an open neighborhood of $c$ in $\CP^1$ and
$O(f(c))$ be an open neighborhood of $f(c)$. Then, by \cite[Theorem
E.1.1]{mas}, there are $C^{2}$-coordinate chart $\phi: (O(c),c)
\rightarrow ( \C^{1}, 0)$ and a $C^{1}$-coordinate chart $\varphi:(
O(f(c)), f(c)) \rightarrow (\C^{2}, 0)$ such that
\begin{itemize}
\item $\varphi_{*}(J(f(c))) = J_{0}$, where $J_{0}$ is a standard
complex structure on $\C^{2}$.
\item $\phi(c) = 0$ and $\varphi(f(c)) = 0$
\item $f: = \varphi \circ f \circ \phi^{-1}: \phi(O(c))
\rightarrow \C^{2}$ is a polynomial in the variable $z$.
\end{itemize}

The \emph{order at c} is a unique integer $k$ such that $\varphi
\circ f \circ \phi^{-1} \in \mathcal{O}_{k} \setminus
\mathcal{O}_{k+1}$. Let $o(c)$ denote [(the order at $c$)-1]. \\

Let $H_{c}$ be a vector space defined as follows:
\[H_{c} := \{t v \mid t \in \C^{1}, \frac{f(z_{i})}{ \mid f(z_{i}) \mid} \rightarrow v \hspace{2mm} \text{as} \hspace{2mm}
z_{i} \rightarrow 0 \} \subset \C^2 \]

$H_{c}$ is called \emph{a tangent cone at c}, which will be denoted
by $\C_{c}$. By considering a normal coordinate chart, we can regard
the tangent cone $\C_{c}$ at $c$ as a subspace of $T_{f(c)}X$.
$\C_{c}$ is preserved by a $J$-action. See \cite[p626]{mas}. The
following Lemma is from \cite[Sec.1.5.]{she}.

\begin{lemma}
Let $c_{i}$ be a critical points of a $J$-holomorphic map $f$. Then,

\[ \mathcal{O}(f^{*}TX)/df(\mathcal{O}(T \CP^1)) \cong N_{f} \oplus \bigoplus_{i}
\C_{c_{i}}^{o(c_{i})},
\]
where $N_{f}$ is a locally free sheaf.

\end{lemma}

$J$ structure on $f^{*}TX$ induces a $J$ structure on $N_{f}$ since
$f$ is a $J$-holomorphic map. The elliptic operator $D_{f,J}$
induces the elliptic operator $D_{f,J}^{N}$ on $N_{f}$. Thus, we get
a two-step elliptic complex:

\begin{align}\label{e;modif}
0 \rightarrow & W^{1,p}(N_{f})\oplus \bigoplus_{i}
\C_{c_{i}}^{o(c_{i})} & \stackrel{\mathbb{D}_{f,J}^{N}}{\rightarrow}
 & L^{p}(\an N_{f})\oplus \bigoplus_{i} ((\an c_{i}) \otimes \C_{c_{i}}^{o(c_{i})}
) & \rightarrow 0\\
 &  \hspace{2.5cm} (v ,w) & \mapsto \hspace{2mm} & (D^{N}_{f,J}(v),0) & \notag
\end{align}

Note that $(\an c_{i}) \otimes \C_{c_{i}}^{o(c_{i})}$ is trivial
because $\an c_{i}$ is trivial. The resulting non-trivial
$D$-cohomology groups are the following:

\begin{gather}
\zdc{\CP^1}{N_{f}\oplus \bigoplus_{i} \C_{c_{i}}^{o(c_{i})}}
  := \text{Ker} D^{N}_{f,J} \oplus \bigoplus_{i} \C_{c_{i}}^{o(c_{i})}  \notag \\
\fdc{\CP^{1}}{N_{f}\oplus \bigoplus_{i} \C_{c_{i}}^{o(c_{i})} } :=
\text{Coker} D^{N}_{f,J} \cong  \fdc{\CP^{1}}{N_{f}} \notag
\end{gather}

\par There is a long exact sequence of $D$-cohomology group induced from the short exact
sequence (\ref{e;short}).

\begin{lemma}\label{l;long}\cite[Proposition
2.4.2]{she} The following sequence is an exact sequence of
$D$-cohomology groups:\\
\begin{gather}
 0 \rightarrow H^{0}(\CP^1, T \CP^1) \rightarrow H^{0}_{D}(\CP^1,
f^{*} TX) \rightarrow H^{0}_{D}(\CP^1, N_{f} \oplus
\bigoplus_{i}\C_{c_{i}}^{o(c_{i})} ) \rightarrow \notag \\
 \rightarrow H^{1}(\CP^1, T
\CP^1) \rightarrow H^{1}_{D}(\CP^1, f^{*} TX) \rightarrow
H^{1}_{D}(\CP^1, N_{f}) \rightarrow 0 \notag
\end{gather}
\end{lemma}

\vspace{2mm}

Proof. The Lemma follows by applying the Snake Lemma and the
definition of the $D$-cohomology groups to the following exact
sequence of two-step elliptic complexes:

\begin{align}
0 \rightarrow &  W^{1,p}(T \CP^1)& \rightarrow & W^{1,p}(
f^{*}TX)&\rightarrow &  W^{1,p}(
N_{f} \oplus \bigoplus_{i}\C_{c_{i}}^{o(c_{i})} ) & \rightarrow 0 \notag \\
 &  \hspace{0.5cm} \text{\footnotesize{$\overline{\partial}$}} \downarrow & & \hspace{0.5cm} \text{\footnotesize{$D_{f,J}$}}
 \downarrow & & \hspace{0.5cm} \text{\footnotesize{$\mathbb{D}_{f,J}^{N}$}} \downarrow &  \notag\\
0 \rightarrow & L^{p}( \an T\CP^1) & \rightarrow & L^{p}( \an f^{*}
TX) & \rightarrow & \hspace{5mm} L^{p}(\an N_{f}) & \rightarrow 0
\notag
\end{align}

\hfill $\Box$

\vspace{2mm}

\subsection{Tangent Space Splitting Theorem}\label{ss;tansplit}

\par Let $\spa{}$
be the parameter space of genus zero, smooth $(C^{\infty})$, simple
$J$-holomorphic maps, representing the homology class $\beta$. Let
$\cam{}{n}$ be the moduli space of genus zero, $n$-pointed, smooth
$J$-holomorphic stable maps, representing the homology class
$\beta$. Let $f: C \rightarrow
 X$ be a stable map from a reducible (arithmetic) genus 0 curve $C$.
 Then, we call $f$ a \emph{simple map} if the map restricted to each
 irreducible component is simple and none of the irreducible
 components have the same image.
We will denote the subset of irreducible, simple stable maps in
$\cam{}{n}$ by $\am{}{n}^{*}$,
 and the subset of any simple stable maps in $\cam{}{n}$
 by $\cam{}{n}^{*}$. $\cam{}{n} \setminus
 \am{}{n}^{*}$ consists of multiple cover stable maps and
 stable maps having reducible domain curves.\\

\begin{proposition} \label{p;tanofquot}
Let $\text{Aut}(\CP^1)$ act on $\spa{}$ by $f \mapsto f \circ
\varphi^{-1}$, where $\varphi \in \text{Aut}(\CP^1)$. Then, the
following holds: \\
\indent i. the quotient space $\spa{}/ \text{Aut}(\CP^1)$ is
diffeomorphic to $\sam{}{0}$. \\
\indent ii. The tangent space at $[(f,\CP^1)] \in \sam{}{0}$ is
$\zdc{\CP^1}{ \mathcal{N}_{f}} \cong \zdc{\CP^1}{N_{f}}\oplus
\bigoplus_{c_{i}}\C_{c_{i}}^{o(c_{i})}$, where $\mathcal{N}_{f}$,
$N_{f}$, $ c_i$ are the normal sheaf, the normal bundle, the
critical point of $f$ respectively.
\end{proposition}

Proof.  i. is obvious and well-known. See \cite[Chap.6]{mas}. Let's
prove ii.\\
Since $\spa{}$ is a smooth manifold, the action of $Aut(\CP^1)$ is
smooth, free and proper. Let $(f,\CP^1)$ be an element in $\spa{}$
and $O(f)$ be an orbit of $(f, \CP^1)$. By \cite[Lemma B.19]{ggk},
$O(f)$ is a smooth, embedded submanifold. \cite[Lemma B.20]{ggk}
implies that the tangent space of the orbit $O(f)$ at $(f, \CP^1)$
is the Lie algebra $H^{0}(\CP^1, T\CP^1)$ of $Aut(\CP^1)$. The
result follows from the short exact sequence of $D$-cohomology
groups (cf.Lemma \ref{l;long}):
\[0 \rightarrow H^{0}(\CP^1, T \CP^1) \rightarrow H^{0}_{D}(\CP^1,
f^{*} TX) \rightarrow H^{0}_{D}(\CP^1, \mathcal{N}_{f}) \rightarrow
0 \]

and the diffeomorphism in i.   \hfill $\Box$\\

\begin{remark}\label{r;immersion}
The tangent space at $[(f,\CP^1)]$ is $\zdc{\CP^1}{N_{f}}$ if and
only if $f$ is an immersion.
\end{remark}

\begin{lemma}\label{l;regsub}
Let $F_{n}^{i}$ be an $i$-th forgetful map. The $i$-th forgetful map
$F_{n}^{i}:\sam{}{n} \rightarrow \sam{}{n-1}$ is a submersion for
any $n \geq 1$.
\end{lemma}

Proof. Let $[(f, \CP^1, a_{1}, \ldots, a_{n})]$ be in $ \sam{}{n}$.
Let $\lambda: (- \epsilon, \epsilon) \rightarrow \sam{}{n-1}$ be a
smooth path such that $\lambda(t) = [( f^{t}, \CP^1, a_{1}^{t},
\ldots, a_{n-1}^{t})]$, $\lambda(0)= [(f, \CP^1, a_{1}, \ldots,
a_{n-1})]$ and $d \lambda(\frac{\partial}{\partial t}\mid_{t=0})
\neq 0$. There is a canonical lifting $\widetilde{\lambda}$ of the
path defined by $\widetilde{\lambda}(t) = [( f^{t}, \CP^1,
a_{1}^{t}, \ldots, a_{n-1}^{t}, a_{n})]$. We have
$\widetilde{\lambda}(0) = [(f,\CP^1, a_{1}, \ldots, a_{n})]$,
 $d \widetilde{\lambda}(\frac{\partial}{\partial t}\mid_{t=0}) \neq 0$
  and $d \lambda(\frac{\partial}{\partial t}\mid_{t=0}) =
  dF_{n}^{n} \circ d \widetilde{\lambda}(\frac{\partial}{\partial t}\mid_{t=0})$.
Other $i$-th forgetful maps can be proved in the same way. So, the Lemma follows. \hfill $\Box$ \\

\begin{proposition}\label{p;tanreg}
 The \hfill tangent \hfill space \hfill at \hfill $[(f,\CP^1, a_{1}, \ldots, a_{n})]$ \hfill in \\
 $\sam{}{n}$ is
$ \zdc{\CP^1}{N_{f}}\oplus \bigoplus_{c_{i}}\C_{c_{i}}^{o(c_{i})}
\oplus \bigoplus_{a_{i}}T_{a_{i}}\CP^1 $, where $c_i$ is a critical
point of $f$.
\end{proposition}

Proof. Assume \hfill the \hfill Proposition \hfill holds \hfill for \hfill $n-1 \geq 0$. \hfill Since \\
$F_{n}^{n}$ is a submersion by Lemma \ref{l;regsub},
$(F_{n}^{n})^{-1}([(f, \CP^1, a_{1}, \ldots, a_{n-1})])$ is smooth
and diffeomorphic to $\CP^1 \setminus \{ a_{1}, \ldots, a_{n-1} \}$
by  $[(f, \CP^1, a_{1}, \ldots, a_{n})] \mapsto a_{n}$.
 Thus, the
tangent space of the fibre $(F_{n}^{n})^{-1}([(f, \CP^1, a_{1},
\ldots,a_{n-1})])$ at\\
 $[(f, \CP^1, a_{1}, \ldots, a_{n})]$ is
isomorphic to $T_{a_{n}} \CP^1$. The Proposition follows by the
induction assumption and Proposition \ref{p;tanofquot}. \hfill $\Box$\\

The type of stable maps in Definition \ref{d;maps} are based on the
types of singularities on the image curve $f(\CP^1)$. \vspace{2mm}

\begin{defn} \label{d;maps}
Let $\mathbf{f}:=(f, \CP^1 , a_{1}, \ldots, a_{n})$ be a simple
pointed stable map. The \emph{node singularity} of the stable map is
the point $p \in f(\CP^1)$ such that $f^{-1}(p)$ consists of two
points $p_{1}$, $p_{2}$ and $df(T_{p_{1}} \CP^1) \cap df(T_{p_{2}}
\CP^1) =  \{0 \}$. We call $\mathbf{f}$ a \emph{nodal stable map} if
the image curve $f(\CP^1)$ has only node singularities. The
\emph{tacnode singularity} of the stable map is the point $p \in
f(\CP^1)$ such that $f^{-1}(p)$ consists of two points $p_{1}$,
$p_{2}$ and $df(T_{p_{1}} \CP^1) = df(T_{p_{2}} \CP^1)$. We call
$\mathbf{f}$ \emph{a tacnode stable map} if the singularities of the
image curve $f(\CP^1)$ has a unique tacnode singularity and all
other singularities in $f(\CP^1)$ are node singularities. The
\emph{triple node singularity} of the stable map is the point $p \in
f(\CP^1)$ such that $f^{-1}(p)$ consists of three points $p_{1}$,
$p_{2}$, $p_{3}$ and $df(T_{p_{i}}\CP^1) \cap df(T_{p_{j}} \CP^1) =
\{ 0 \}$ if $i \neq j$. We call $\mathbf{f}$ a \emph{triple node
stable map} if $f$ has a unique triple node singularity and all
other singularities are node singularities. The \emph{cuspidal
singularity } $p \in f(\CP^1)$ is the image of an order 2
singularity of $f$. We call the stable map $\mathbf{f}$ a
\emph{cuspidal stable map} if $f(\CP^1)$ has a unique cuspidal
singularity and all other singularities are node singularities. We
will call the two stable $J$-holomorphic maps $f, g$ are
\emph{equi-singular} if they have the same number of the irreducible
components on the domain curve and they have the same type of
singularities in their image curves.
\end{defn}

\vspace{2mm}

Readers may refer to \cite{kwona} for the local figures of
singularities in $f(\CP^1)$.\\

\begin{proposition}\label{p;nodalopen}
The nodal stable maps locus is open in $\cam{}{n}$.
\end{proposition}

\vspace{2mm}

Proof. The transversality is an open condition. Thus, a small
perturbation of a nodal stable map doesn't change the type of
singularities in the image curve. A small perturbation doesn't
create other type of singularities because of the Symplectic
adjunction formula and the type of singularities in the image curve
of a nodal stable map.
The result is obvious. \hfill $\Box$ \\

\begin{proposition}\label{p;cuspcodim}~\cite[Theorem 3.2.1]{she}.
The cuspidal stable map's locus is a real codimension two locus in
$\cam{}{n}$.
\end{proposition}

\begin{remark}\label{r;cusptangent}(Deformation property near the
degree 2 singularity) Let $f$ be a cuspidal stable map with a
critical point at $c$. Then, the tangent space at $\mathbf{f} :=
[(f,\CP^1)]$ is $\zdc{\CP^1}{N_{f}}\oplus \C_{c} $. After the
coordinate changes on the neighborhood of $c$ in the domain curve
and the neighborhood of $f(c)$ in the target space, $f$ is
represented by $z \mapsto (z^{2}, z^{3})$. Since
$\zdc{\CP^{1}}{N_{f}}$ parameterizes a first order deformation, any
smooth path $\xi: (-\epsilon, \epsilon) \rightarrow \sam{}{n}$
tangential to a vector in $\zdc{\CP^{1}}{N_{f}}$ and $\xi(0) =
\mathbf{f}$ parameterizes a first order deformation $f_{\eta}(z)=
(z^{2}, z^{3}+ \eta z)$ for $\eta \in (-\epsilon, \epsilon)$ after
the coordinate changes. One can check that $f_{\eta}$ has a node
singularity if $\eta \neq 0$. Thus, Proposition \ref{p;cuspcodim}
shows that $\C_{c}$ is tangential to the cuspidal stable maps locus.
By the definition of the cuspidal stable map, $\C_{c}$ parameterizes
a \emph{second order deformation} near the singularity $c$.
\end{remark}

\begin{lemma}\label{l;difim}
Let $\mathbf{f} := [(f,\CP^1, a_{1}, \ldots, a_{n})]$ be in
$\sam{}{n}$. Let $ev_{i}$ denote the $i$-th evaluation map, defined
by $ev_{i}(\mathbf{f}) = f(a_{i})$. Then, the differential of the
$i$-th evaluation map
 at $\mathbf{f}$ is given as follows:
\begin{gather}
[\zdc{\CP^1}{N_{f}}\oplus  \bigoplus_{c_{j}} \C_{c_{j}} ] \oplus
\bigoplus_{a_{j}}T_{a_{j}}\CP^1 \rightarrow T_{f(a_{i})}X  \notag \\
(v_{1} + v_{2}) + w  \mapsto v_{1}(a_{i}) +
df_{a_{i}}(\pi_{i}(w))\label{e;differential}
\end{gather}
where $v_{1}$, $v_{2}$, $w$ belong to $\zdc{\CP^1}{N_{f}}$,
$\bigoplus_{c_{j}}\C_{c_{j}}$, $\bigoplus_{a_{j}}T_{a_{j}}\CP^1$
respectively, and $\pi_{i}:\bigoplus_{a_{j}}T_{a_{j}} \rightarrow
T_{a_{i}} \CP^{1}$ is a natural projection map to the $i$-th
component.
\end{lemma}

Proof. Since $v_{1}$, $v_{2}$, $w$ are independent vectors, it is
enough to check $v_{1} \mapsto v_{1}(a_{i})$, $v_{2} \mapsto 0 $, $w
\mapsto df_{a_{i}}(w) $. Let $\gamma:(-\epsilon, \epsilon)
\rightarrow \sam{}{n}$ be a smooth path tangential to $v_{1}$ and
$\gamma(0) = \mathbf{f}$. $ d(ev_{i} \circ
\gamma)(\frac{\partial}{\partial t} \mid_{t=0}) =$$\displaystyle
\frac{d}{dt}exp(t \cdot v_{1}(a_{i} )) = v_{1}(a_{i}) $.
$\C_{c_{j}}$ parameterizes the 2nd order deformation. Thus, $v_{2}
\mapsto 0$. It is obvious that if $w \in T_{a_{j}} \CP^1$, $j \neq
i$, then $dev_{i}(w) =0$ because $T_{a_{j}} \CP^1$ parameterizes the
deformation of the $j$-th marked point while we consider the $i$-th
evaluation map. Let $\lambda:(-\varepsilon,\varepsilon) \rightarrow
\CP^1$ be a smooth path such that $\lambda(0)= a_{i}$, $d
\lambda(\frac{\partial}{\partial t} \mid_{t=0}) = w \in T_{a_{i}}
\CP^1$. Then, the induced path $ev_{i} \circ \lambda: (-\varepsilon
,\varepsilon) \rightarrow X$ satisfies that $ev_{i} \circ \lambda(0)
= f(a_{i})$ and $d(ev_{i} \circ \lambda)(\frac{\partial}{\partial t}
\mid_{t=0}) = df_{a_{i}}(\pi_{i}(w))$.
 \hfill $\Box$ \\

\vspace{2mm}

\par Let $C$ be a pointed reducible curve which has two irreducible
components $(C_{1}, a_{1}, \ldots, a_{r})$, $(C_{2}, b_{1}, \ldots,
b_{s})$, where $r + s = n$. Let $q_{1} \in C_{1}$, $q_{2} \in C_{2}$
be pregluing  points.
 Then, a pointed stable map
$(f,C,a_{1}, \ldots, a_{r}, b_{1}, \ldots, b_{s})$ can be  written
as $((f_{1}, (C_{1}, q_{1}), a_{1}, \ldots, a_{r}), (f_{2}, (C_{2},
q_{2}), b_{1},\ldots, b_{s}))$.  \\
Let $\red$ denote the set of stable maps in $\cam{}{n}$ such that
\begin{itemize}
\item $[((f_{1}, (C_{1}, q_{1}), a_{1},
\ldots, a_{r}), (f_{2}, (C_{2}, q_{2}), b_{1},\ldots, b_{s}))]$,
where $C_{i} \cong \CP^1$
\item $f_{*}([C_{1}])= \beta_{1}, \hspace{0.5cm} f_{*}([C_{2}]) =
\beta_{2}$
\end{itemize}
$\red$ is generically smooth and forms a (real) codimension two
subspace in $\cam{}{n}$.

\begin{proposition}\label{p;2redtan}
Let's consider a stable map \\
$f_{1} \bigvee f_{2} := [((f_{1}, (C_{1}, q_{1}), a_{1}, \ldots,
a_{r}), (f_{2}, (C_{2},
q_{2}), b_{1}, \ldots, b_{s}))] \in \red$.\\

\noindent (i) Let's assume that $f_{i}$, $i=1,2$, is an immersion.
 Then,
the tangent space splitting at $f_{1} \bigvee f_{2}$ is: \\

\begin{gather}
\zdc{C_{1}}{ N_{1}} \oplus \zdc{C_{2}}{N_{2}} \oplus T_{a_{1}}C_{1}
\oplus \ldots \oplus T_{a_{r}}C_{1} \oplus T_{b_{1}}C_{2} \oplus
\ldots \oplus T_{b_{s}}C_{2} \oplus \notag \\
\oplus (T_{q_{1}}C_{1} \otimes T_{q_{2}}C_{2}) \oplus T_{q_{1}}C_{1}
\oplus T_{q_{2}}C_{2} \ominus T_{f(q)} X \notag
\end{gather}
\vspace{2mm}

\noindent where $N_{i}$ is a normal bundle and $q$ is a node in $C$.\\

\noindent (ii) Let's assume $f_{1*}([C_{1}])$ is trivial and $f_{2}$
is an immersion.
 Then, the tangent space splitting at $f_{1} \bigvee f_{2}$ is:\\
\begin{gather}
 H^{1}(C_{1}, T C_{1}(-q_{1} -a_{1} - \ldots -
a_{r} )) \oplus \zdc{C_{2}}{ N_{2}} \oplus \notag \\
\oplus T_{b_{1}}C_{2} \oplus \ldots \oplus T_{b_{s}}C_{2} \oplus
(T_{q_{1}}C_{1} \otimes T_{q_{2}}C_{2})\oplus T_{q_{2}}C_{2} \notag
\end{gather}
\end{proposition}

Proof. (i)Let $M_{r+1}(X , \beta_{1}, J)^{**}$, $M_{s+1}(X ,
\beta_{2}, J)^{**}$ be the open subset of $M_{r+1}(X, \beta_{1},J)$,
$M_{s+1}(X, \beta_{2},J)$ consisting of stable maps without critical
points. Let's consider a smooth map $ev_{\beta_{1}} \times
ev_{\beta_{2}}$:
\begin{gather}
 M_{r+1}(X , \beta_{1}, J)^{**} \times
M_{s+1}(X , \beta_{2}, J)^{**} \rightarrow X \times X \notag \\
([(f_{1},(C_{1},q_{1}),a_{1},
\ldots,a_{r})],[(f_{2},(C_{2},q_{2}),b_{1}, \ldots,b_{s})]) \mapsto
(f_{1}(q_{1}), f_{2}(q_{2})) \notag
\end{gather}
Let $(\red)^{**}$ be the subset of $\red$ consisting of the stable
maps whose restriction to each irreducible component is an
immersion. Then, $(\red)^{**}$ is diffeomorphic to
$ev_{\beta_{1}}^{-1} \times ev_{\beta_{2}}^{-1}(\{(q,q) \mid q \in X
\})$. Let's denote $[(f_{1},(C_{1},q_{1}),a_{1}, \ldots,a_{r})]$,
$[(f_{2},(C_{2},q_{2}),b_{1}, \ldots,b_{s})]$ by $\mathbf{f}_{1}$,
$\mathbf{f}_{2}$ respectively. By Proposition \ref{p;tanreg}, we
have:
\begin{gather}
T_{\mathbf{f}_{1}}M_{r+1}(X, \beta_{1}, J)^{**}
 \cong \zdc{C_{1}}{N_{1}} \oplus T_{q_{1}}C_{1} \oplus \bigoplus_{a_{i}} T_{a_{i}}C_{1} \label{e;tanq1} \\
T_{\mathbf{f}_{2}}M_{s+1}(X, \beta_{2}, J)^{**} \cong
\zdc{C_{1}}{N_{2}} \oplus T_{q_{2}}C_{2} \oplus \bigoplus_{b_{i}}
T_{b_{i}}C_{2} \label{e;tanq2}
\end{gather}

The differential $dev_{\beta_{1}} \times dev_{\beta_{2}}$ at
$(\mathbf{f}_{1},\mathbf{f}_{2})$ is as follows:

\begin{gather} T_{\mathbf{f}_{1}}M_{r+1}(X, \beta_{1}, J)^{**} \times
T_{\mathbf{f}_{2}}M_{s+1}(X, \beta_{2}, J)^{**}
\rightarrow  T_{f(q_{1})}X \times T_{f(q_{2})}X  \notag \\
(v + \xi + \varsigma , v' + \xi'+ \varsigma')  \mapsto  (v(q_{1}) +
df_{1q_{1}}(\xi), v'(q_{2}) + df_{2q_{2}}(\xi')) \notag
\end{gather}
where $v, v'$ are elements in $\zdc{C_{1}}{N_{f_{1}}}$,
$\zdc{C_{2}}{N_{f_{2}}}$ respectively, $\xi, \xi'$ are elements in
$T_{q_{1}}C_{1}$, $T_{q_{2}}C_{2}$ respectively, $\varsigma$,
$\varsigma'$ are elements in $\bigoplus_{a_{i}} T_{a_{i}}C_{1}$,
$\bigoplus_{b_{i}} T_{b_{i}}C_{2}$ respectively. By Lemma
\ref{l;difim}, $ev_{\beta_{i}}$ is a submersion. Thus, the
differential $dev_{\beta_{1}} \times dev_{\beta_{2}}$ at any
$(\mathbf{f}_{1},\mathbf{f}_{2})$ is surjective. Let
$(\mathbf{f}_{1},\mathbf{f}_{2}) = f_{1} \bigvee f_{2}$, i.e.,
$f_{1}(q_{1}) = f_{2}(q_{2})$. By combing with the natural
surjective map $T_{f(q)}X \times T_{f(q)}X \rightarrow T_{f(q)}X$,
$(w, w') \mapsto w - w'$, we get the following surjective map:
\begin{gather}
 T_{\mathbf{f}_{1}}M_{r+1}(X, \beta_{1}, J)^{**} \times
T_{\mathbf{f}_{2}}M_{s+1}(X, \beta_{2}, J)^{**} \rightarrow
T_{f(q)}X \label{e;comp}\\
(v + \xi + \varsigma , v' + \xi'+ \varsigma')  \mapsto v(q_{1}) +
df_{1q_{1}}(\xi)- v'(q_{2}) - df_{2q_{2}}(\xi') \notag
\end{gather}

By Implicit Function Theorem, $(\red)^{**}$ is smooth manifold of
dimension $\text{dim}M_{r+1}(X, \beta_{1}, J)^{**} +
\text{dim}M_{s+1}(X, \beta_{2}, J)^{**} - \text{dim}X$ and the
kernel of the map (\ref{e;comp}) is the tangent space of
$(\red)^{**}$ at $(f_{1} \bigvee f_{2})^{**}$. The statement (i) is
from a $K$-group expression of the orthogonal decomposition induced
by the following short exact sequence (\ref{e;sho}), combining with
(\ref{e;tanq1}),(\ref{e;tanq2}). The term $T_{q_{1}}C_{1} \otimes
T_{q_{2}}C_{2}$ parameterizes a smoothing node deformation.

\begin{multline} \label{e;sho}
0 \rightarrow T_{f_{1} \bigvee f_{2}}(\red)^{**} \rightarrow \\
\rightarrow T_{\mathbf{f}_{1}}M_{r+1}(X, \beta_{1}, J)^{**} \times
T_{\mathbf{f}_{2}}M_{s+1}(X, \beta_{2}, J)^{**} \rightarrow
T_{f(q)}X \rightarrow 0
\end{multline}

\par (ii) Since $f_{1*}([C_{1}])$ is trivial and a stable map,
the irreducible component $C_{1}$ has at least 3 marked points. The
map $M_{r+1}(X, 0, J) \rightarrow M_{r+1} \times X$, $\mathbf{f}_{1}
\mapsto ([(C_{1}, q_{1}, a_{1}, \ldots, a_{r})], f_{1}(q_{1}))$ is a
diffeomorphism, where $M_{r+1}$ is the Deligne-Mumford moduli space
of $r+1$-pointed smooth genus 0 curves. The tangent space at
$[(C_{1}, q_{1}, a_{1}, \ldots, a_{r})]$ in $M_{r+1}$ is
$H^{1}(C_{1}, TC_{1}(-q_{1}-a_{1} - \ldots - a_{r}))$ which is the
space of first order deformations of a pointed smooth Riemann
surface. Thus, the tangent space at $\mathbf{f}_{1}$ is:
\begin{equation}\label{e;tantri}
T_{\mathbf{f}_{1}}M_{r+1}(X, 0, J) \cong H^{1}(C_{1},
TC_{1}(-q_{1}-a_{1} - \ldots - a_{r})) \oplus T_{f(q)}X
\end{equation}

With the same notation we used in (i) for $ v', \xi', \varsigma'$
and the similar argument, we get a submersion map:
\begin{gather}
 T_{\mathbf{f}_{1}}M_{r+1}(X, \beta_{1}, J) \times
T_{\mathbf{f}_{2}}M_{s+1}(X, \beta_{2}, J)^{**} \rightarrow
T_{f(q)}X \label{e;compt}\\
(\varsigma + \xi , v' + \xi'+ \varsigma')  \mapsto \xi - v'(q_{2}) -
df_{2q_{2}}(\xi') \notag
\end{gather}

where $\varsigma \in H^{1}(C_{1}, TC_{1}(-q_{1}-a_{1} - \ldots -
a_{r}))$, $\xi \in T_{f(q)}X$. The rest of arguments are the same
with (i). \hfill $\Box$

\vspace{2mm}

Repeated similar calculations result in the following Theorem. See
\cite{kwont} for the calculation in algebraic geometry.

\begin{theorem}\label{t;fulltangsp}\emph{(Tangent Space Splitting
Theorem)} \\
Let $[(f,C, a_{1}, \ldots, a_{n})]$ be a point in $\cam{}{n}$ such
that $f$ is an immersion on each irreducible component. Let
$\tilde{p}:\widetilde{C}:= \bigsqcup_{i=1}^{r}C_{i} \rightarrow C$
be a normalization map of $C$. Then, the tangent space splitting at
$[(f,C, a_{1}, \ldots, a_{n})]$ is:
\begin{gather}
\bigoplus_{i=1}^{r}\zdc{C_{i}}{N_{i}} \oplus
\bigoplus_{i=1}^{n}T_{a_{i}}C_{i} \oplus
[\bigoplus_{i=1}^{\delta}(T_{p_{i}}C_{\nu(p_{i})} \otimes
T_{p'_{i}}C_{\nu(p'_{i})})] \oplus \notag \\
\oplus [\bigoplus_{i=1}^{\delta}T_{p_{i}}C_{\nu(p_{i})} \oplus
\bigoplus_{i=1}^{\delta}T_{p'_{i}}C_{\nu(p'_{i})}] \ominus
[\bigoplus_{i=1}^{\delta}T_{f(q_{i})}X]  \label{e;genetan}
\end{gather}
where $\delta$ is the number of gluing points, $p_{i} \in
C_{\nu(p_{i})} $, $p'_{i} \in C_{\nu(p'_{i})}$, $i=1, \ldots,
\delta$, $\{\nu(p_{i}), \nu(p'_{i}) \} \subset \{1, \ldots, r \}$
are pregluing points such that $\tilde{p}(p_{i})=
\tilde{p}(p'_{i})$, and $q_{i}$, $i=1, \ldots, \delta$, is a gluing
point.
\end{theorem}

\section{Singularity Analysis and Calculations of Ramification Indices} \label{s;index}

\begin{lemma}\label{l;kerofdev} Let $\mathbf{f} := [(f,\CP^{1},a_{1}, \ldots,  a_{n})] \in \cam{}{n}$
be a point represented by a stable map with $l$($\geq 0$) singular
points of order 2. The kernel of the differential $dev := d ( ev_{1}
\times \ldots \times ev_{n})$ at $\mathbf{f}$ is isomorphic to
$\zdc{\CP^1}{N_{f}(-a_{1}- \ldots - a_{n})} \oplus
\bigoplus_{j=1}^{l} \C_{c_{j}}$, where $a_{i}$, $i=1, \ldots, n$, is
a regular point of $f$.
\end{lemma}

Proof. By Lemma \ref{l;difim}, it is clear that $\bigoplus_{j=1}^{l}
\C_{c_{j}}$ is in the kernel of $dev$. $v_{1}(a_{i})$, $df(
\pi_{i}(w))$ in (\ref{e;differential}) are independent vectors.
Thus, $dev(v_{1} + v_{2} + w) = 0 $ if and only if $v_{1}(a_{i}) =
0$ and $df( \pi_{i}(w))=0$ for all $i=1, \ldots, n$.\\
If $v_{1}(a_{i}) = 0$ for $i=1, \ldots, n$, then $v_{1} \in
\zdc{\CP^1}{N_{f}(-a_{1}- \ldots - a_{n})}$ and vice versa. Since
$a_{i}$ is a non-singular point of $f$, $\pi_{i}(w)$ is a zero
vector iff $df(\pi_{i}(w))$ vanishes. Thus, the result follows.
\hfill $\Box$\\

\begin{proposition}\label{p;regularity}
Let $n := c_{1}(f^{*}TX) - 1$. Let $\mathbf{f}:= [(f,\CP^1, a_{1},
\ldots, a_{n})]$ be a point in $\cam{}{n}$ such that $f$ is an
immersion. Then, the map $ev:= ev_{1} \times \ldots \times ev_{n}$
is regular at $\mathbf{f}$.
\end{proposition}

Proof. Simple dimension count shows dim $\zdc{\CP^1}{N_{f}(-a_{1}-
\ldots - a_{n})}= 0$. Therefore, $dev$ is one-to-one at
$\mathbf{f}$.
 Since $T_{\mathbf{f}} \cam{}{n}$ is a finite dimensional
vector space, the result follows from Lemma \ref{l;kerofdev}. \hfill
$\Box$\\

Nodal stable map, triple node stable map, and tac node stable map
are immersions. So, we get the following Corollary.

\begin{corollary}\label{c;exregularity}
Let $n :=  c_{1}(f^{*}TX) - 1$. Let $\mathbf{f}:= [(f,\CP^1, a_{1},
\ldots, a_{n})]$ be a point in $\cam{}{n}$. If $f$ is a nodal stable
map, triple node stable map, or tac node stable map, then the map
$ev:= ev_{1} \times \ldots \times ev_{n}$ is regular at
$\mathbf{f}$.
\end{corollary}

\begin{lemma}\label{l;cokofdev} Let $n := c_{1}(f^{*}TX) - 1$.
The cokernel of $dev$ at $\mathbf{f} := [(f,\CP^{1},a_{1},
\ldots,a_{n})] $ $\in \cam{}{n}$ is isomorphic to
 $\fdc{\CP^1}{N_{f}(-a_{1}- \ldots - a_{n})}$,
where $N_{f}$ is the normal bundle and $a_{i}$, $i=1, \ldots, n$, is
a non-singular point of $f$.
\end{lemma}

Proof. By Lemma \ref{l;difim}, we have

\[ \mbox{coker} dev  \cong  \frac{ \bigoplus_{i=1}^{n} T_{f(a_{i})}
X }{dev( \zdc{ \CP^1}{ N_{f}}) \oplus dev( \bigoplus_{i=1}^{n}
T_{a_{i}} \CP^1)}\]

Let $B$ and $C$ be subvector spaces of the vector space $A$ such
that $B \cap C = \{ 0 \}$. Then, one can easily check the elementary
isomorphism $ \displaystyle \frac{A}{ B \oplus C } \cong
\frac{C^{\perp} \oplus C }{B \oplus C} \cong \frac{C^{\perp}}{B} $,
where $C^{\perp}$ is an orthogonal complement of $C$ in $A$.

Consider the short exact sequence of sheaves

\begin{equation}\label{e;shortsheaf}
0 \rightarrow \mathcal{O}(N_{f}(-a_{1}- \ldots - a_{n})) \rightarrow
\mathcal{O}(N_{f}) \rightarrow \bigoplus_{i=1}^{n} \nu_{i}
\rightarrow 0,
\end{equation}
 where $\nu_{i}$ is a skyscraper sheaf
supported by $a_{i}$, $i=1, \ldots, n$.

Since $dev( \zdc{ \CP^1}{ N_{f}}) \cap dev( \bigoplus_{i=1}^{n}
T_{a_{i}} \CP^1) \cong \{ 0 \} $, we get

\begin{align}
 \mbox{coker} dev & \cong (\bigoplus_{i=1}^{n} N_{a_{i}}) / dev(\zdc{ \CP^1}{ N_{f}})\\
& \cong \zdc{\CP^1}{\bigoplus_{i=1}^{n} \nu_{i}}/ dev(\zdc{ \CP^1}{
N_{f}}), \label{e;qu}
\end{align}

By Lemma \ref{l;long}, (\ref{e;shortsheaf}) induces a long exact
sequence of D-cohomology groups:

\begin{gather*}
0 \rightarrow \zdc{\CP^1}{ N_{f}(-a_{1}- \ldots - a_{n})}
\rightarrow \zdc{\CP^1}{ N_{f}} \rightarrow \\
\rightarrow \zdc{\CP^1}{ \bigoplus_{i=1}^{n} \nu_{i} } \rightarrow
\fdc{ \CP^1}{ N_{f}(-a_{1}- \ldots - a_{n})} \rightarrow \\
\rightarrow \fdc{\CP^{1}}{ N_{f}}  \rightarrow \ldots.
\end{gather*}

By Riemann-Roch's Theorem, $\zdc{\CP^1}{ N_{f}(-a_{1}- \ldots -
a_{n})}$ and $\fdc{\CP^{1}}{ N_{f}}$ vanish. The Lemma follows.
\hfill $\Box$

\begin{theorem}\label{t;cuspind}
Let $n := c_{1}(f^{*}TX) - 1$. Let $\mathbf{f} := [(f,\CP^1,a_{1},
\ldots, a_{n})]$ be represented by a cuspidal stable map in
$\cam{}{n}$, where $a_{i}$ is a regular point. Then, $\mathbf{f}$ is
a critical point of the $ev$ map. The ramification index of the ev
map along the cuspidal stable maps locus is 2.
\end{theorem}

Proof. Lemma \ref{l;kerofdev} implies the kernel of
$dev_{\textbf{c}} : = \tau$, where $\tau$ is a skyscraper sheaf
supported by the cuspidal singularity. The cokernel of $dev$ is \\
$\zdc{\CP^1}{N_{f}(-a_{1}- \ldots - a_{n})}$ by Lemma
\ref{l;cokofdev}. The Kodaira-Serre duality in \cite[Lemma
1.5.1]{she} shows that $\fdc{\CP^1}{ N_{f}(-a_{1}- \ldots - a_{n})}$
is isomorphic to $\zdc{\CP^1}{ \omega_{\CP^1} \otimes [N_{f}(-a_{1}-
\ldots - a_{n})]^{*}}$, where $\omega_{\CP^1}$ is a dualizing sheaf
on $\CP^1$ and $[N_{f}(-a_{1}- \ldots - a_{n})]^{*}$ is a dual
vector bundle of $[N_{f}(-a_{1}- \ldots - a_{n})]$. There is a
canonical residue morphism of degree 2 from the local slice of the
direction $\tau$ to the direction  $\fdc{\CP^1}{ N_{f}(-a_{1}-
\ldots - a_{n})}$ induced by the $ev$ map. By Micallef-White's
Theorem, after the coordinate chages, the morphism can be written:

\begin{align}
 \tau \simeq \tau^{*} &
\rightarrow & \zdc{\CP^1}{ \omega_{\CP^1} \otimes [N_{f}(-a_{1}-
\ldots -
a_{n})]^{*}}^{*} \notag \\
v dz & \mapsto & \frac{1}{2 \pi i}\int_{\gamma} \frac{v^{2}}{z} dz
\hspace{5cm}  \notag
\end{align}
because the local index of $f$ is 2 and $\tau$ parameterizes the 2nd
order deformation. Thus, the Theorem follows.
 \hfill $\Box$ \\

\begin{proposition}\label{p;red}
Let $n := c_{1}(f^{*}TX) - 1$ and $n = r + s$. Let the point
\begin{align}
\mathbf{f}& :=[(f,C,z_{1},\ldots, z_{r}, w_{1}, \ldots, w_{s})] \notag \\
& := [((f_{1}, (C_{1}, q_{1}), z_{1}, \ldots, z_{r}), (f_{2},
(C_{2}, q_{2}), w_{1}, \ldots, w_{s}))] \notag
\end{align}
be represented by a reducible stable map $f$, where $f_{i}$ is an
immersion if $f_{i}$ is not trivial. Then, \\
$(i)$ If any of $f_{i}$ is a degree 0 map, then the cokernel of
$dev$ at $\mathbf{f}$ has a (real) rank at least six. \\
$(ii)$ If $r$ or $s$ is strictly bigger than $c_{1}(f_{1}^{*}TX) -2$
or $c_{1}(f_{2}^{*}TX) -2$ respectively, then the cokernel of $dev$
at $\mathbf{f}$ has a (real) rank bigger
than four. \\
$(iii)$ If $r$ or $s$ is $c_{1}(f_{1}^{*}TX)-2$ or
$c_{1}(f_{2}^{*}TX) -2$ respectively,
then the cokernel of $dev$ at $\mathbf{f}$ has a (real) rank two.\\
$(iv)$ If $r$ or $s$ is $c_{1}(f_{1}^{*}TX) -1$ or
$c_{1}(f_{2}^{*}TX)-1$ respectively, then the evaluation map $ev$ at
$\mathbf{f}$ is regular.
\end{proposition}
\vspace{2mm}

Sketch of the Proof. In Proposition \ref{p;2redtan} (i), the vector
space $ T_{q_{1}}C_{1} \oplus T_{q_{2}}C_{2} \ominus T_{f(q)} X$
parameterizes the node deformation. So, it doesn't contribute to the
rank of the ev map. $T_{q_{1}}C_{1} \otimes T_{q_{2}}C_{2}$
generates the first order deformation of smoothing node. During the
smoothing node deformation, the stable map also changes because the
stable maps are the same if they agree to infinite order at any
point. Since the deformation space generating smoothing node
deformation is smooth and parameterize the 1st order deformation,
the image of any non-trivial vector in $T_{q_{1}}C_{1} \otimes
T_{q_{2}}C_{2}$ by $dev$ is non-trivial. Therefore, the (real) rank
dimension contribution from $T_{q_{1}}C_{1} \otimes T_{q_{2}}C_{2}$
for the ev map is two. Let's prove (i). Suppose that $f_{1}$ is a
trivial map. The stability condition implies that $C_{1}$ contains
at least 2 marked points. The differential of the $i$-th evaluation
map $ev_{i}$ is zero on $C_{1}$. The maximum dimensional
contribution to the rank of the $ev$ map from $C_{2}$ is at most
$2[(c_{1}(f^{*}TX)-3) + (c_{1}(f^{*}TX)-3)]$. The first
$2(c_{1}(f^{*}TX)-3)$ is contributed by $T_{w_{1}} C_{2}, \ldots,
T_{w_{s}}C_{2}$, where $s = c_{1}(f^{*}TX)-3$. The second
$2(c_{1}(f^{*}TX)-3)$ is contributed by $\zdc{C_{2}}{N_{f_{2}}}$.
The contribution toward the rank of the ev map from
$\zdc{C_{2}}{N_{f_{2}}}$ is limited by the number of marked points
$s$ in $C_{2}$. The comparison with the dimension
$4(c_{1}(f^{*}TX)-1)$ of the moduli space leads to the result (i).
The proofs of (ii), (iii), (iv)are similar and straightforward. We
will omit them. \hfill $\Box$

\section{Transversality Properties on $\cam{}{n}$}\label{s;transversalityofpoint}

\begin{lemma}\label{l;surj}
Let $n := c_{1}(f^{*}TX) -1$. Let $\mathbf{f} := [(f,C, a_{1},
\ldots, a_{n})]$ be in $\cam{}{n}$. Assume that $f_{j}$ is regular
at $a_{i}$ in the irreducible component $C_{j}$ in $C$. Let $c_{k}$,
$k=1, \ldots, m$, be critical points in $C_{j}$. If $\sum_{k=1}^{m}
o(c_{k})$ (cf. def. of $o(c_{k})$ in sec.\ref{ss;normalchom}) is
less than $ c_{1}(f_{j}^{*}TX) - 1$, then the $i$-th evaluation map
$ev_{i}$ is regular at $\mathbf{f}$.
\end{lemma}

\vspace{2mm}

Proof. The condition in the sum of the order of the critical points
and the dimension count by the Riemann-Roch's Theorem with
Proposition \ref{p;tanreg} implies $\zdc{C_{j}}{N_{f_{j}}}$ has
dimension at least two. Both $dev_{i}(T_{f(a_{i})}X)$ and
$dev_{i}(\zdc{C_{j}}{N_{f_{j}}})$ have dimension two because $f$ is
regular at $a_{i}$. Moreover, $dev_{i}(T_{f(a_{i})}X)$ and
$dev_{i}(\zdc{C_{j}}{N_{f_{j}}})$ are independent vector spaces.
Thus, $ev_{i}$ is regular at $\mathbf{f}$. \hfill $\Box$

\vspace{2mm}

\par Let \hfill $O(\mathbf{f})$ \hfill be \hfill a \hfill small \hfill open
\hfill neighborhood \hfill  of \hfill $\mathbf{f} \in
\bigcap_{i=1}^{n}
ev_{i}^{-1}(q_{i})$ \hfill such \hfill that \\
$O(\mathbf{f}) \bigcap (\bigcap_{i=1}^{n} ev_{i}^{-1}(q_{i})) = \{
 \mathbf{f} \}$. Let $O(q_{i})$ be a small open neighborhood of $q_{i} \in X$. Then,

\[ \mbox{Intersection multiplicity
at} \hspace{2mm} \mathbf{f}:= \#  [O(\mathbf{f}) \bigcap
(\bigcap_{i=1}^{n} ev_{i}^{-1}(q'_{i}))],     \] where $q_{i}'$ is a
generic point in $O(q_{i})$. Lemma \ref{l;surj} shows that the
pull-back cycle $ev_{i}^{-1}(q_{i})$ is smooth on general points of
$ev_{i}^{-1}(q_{i})$. Since the transversality property is an open
condition, the intersection multiplicity is well-defined. We say
\emph{the cycles $ev_{1}^{-1}(q_{1}), \ldots, ev_{n}^{-1}(q_{n})$
meet transversally} if the intersection multiplicity at any point in
$\bigcap_{i=1}^{n} ev_{i}^{-1}(q_{i})$ is one.

\begin{lemma}\label{t;maintransversality}
Let $n := c_{1}(f^{*}TX) -1$. Let $\mathbf{f}$ be in
$\bigcap_{i=1}^{n} ev_{i}^{-1}(q_{i})$, where $q_{i}$,
$i=1, \ldots, n$, are points in general position in the compact semipositive symplectic 4-manifold $X$. Then, the following holds.\\
(i) If $\mathbf{f}$ is represented by a stable map which is an
immersion and has an irreducible domain curve,
 then the intersection multiplicity at $\mathbf{f}$ is one.\\
(ii)If $\mathbf{f}$ is represented by a cuspidal stable map none of
whose marked points is a critical point, then the intersection
multiplicity at $\mathbf{f}$ is two.
\end{lemma}

Proof. Note that $ev^{-1}(q_{1}, \ldots, q_{n}) = \bigcap_{i=1}^{n}
ev_{i}^{-1}(q_{i})$. (i) follows from Proposition \ref{p;regularity}
because $ev$ at $\mathbf{f}$ is local diffeomorphism. (ii) follows
from Theorem \ref{t;cuspind} because the local ramification index is
identical to the
multiplicity.  \hfill $\Box$\\

\vspace{2mm}

All of the results so far sum up to the non-transversality property
in Theorem \ref{t;cusmul} and prove the Tian's deep conjecture.

\begin{theorem}\label{t;cusmul} Let $n := c_{1}(f^{*}TX) -1$. The cuspidal stable maps locus is
the unique equi-singular locus in $\am{}{n}$ of real codimension
$\leq 2$ on which transversality uniformly fails.
\end{theorem}

Proof. The Theorem follows immediately from Lemma
\ref{t;maintransversality}. \hfill $\Box$\\

\vspace{5mm}

\textbf{Acknowledgement} Gang Tian suggested me to write the
symplectic version of the paper \cite{kwona}. Fundamental
contribution of Prof. Tian was relating the intersection theoretic
property to the existence of critical points of the stable map. I
appreciate his explanations about his observations at MIT in June
2003 during his invitation. I thank Selman Akbulut, Yael Karshon for
helpful conversations, e-mail correspondences. I appreciate the nice
working circumstances at the Fields Institute.

\vspace{1cm}

\begin{center}
University of Montana - Western \\
Department of Mathematics \\
710 South Atlantic Street \\
Dillon, MT 59725, USA
\end{center}

\begin{center}
 s$_{-}$kwon@hotdawg.umwestern.edu
\end{center}

\end{document}